\newif\ifpictures
\picturestrue 

\documentclass[12pt]{amsart}
\usepackage{amssymb}
\usepackage{epsf}
\usepackage[dvips]{graphicx}

\numberwithin{equation}{section}

\newtheorem{thm}{Theorem}
\newtheorem{prop}[thm]{Proposition}
\newtheorem{lemma}[thm]{Lemma}
\newtheorem{cor}[thm]{Corollary}

\newtheorem{example}[thm]{Example}
\newtheorem{remark}[thm]{Remark}

\newenvironment{rem}{\begin{remark}\rm}{\end{remark}}
\newcounter{FNC}[page]
\def\newfootnote#1{{\addtocounter{FNC}{2}$^\fnsymbol{FNC}$%
     \let\thefootnote\relax\footnotetext{$^\fnsymbol{FNC}$#1}}}

\newcommand{\R}{\mathbb{R}}
\newcommand{\E}{\mathbb{E}}

\newcommand{\bv}{{v}}

\newcommand{\stirlingnumbersecond}[2]{\genfrac{\{}{\}}{0pt}{}{#1}{#2}}
\newcommand{\atopfrac}[2]{\genfrac{}{}{0pt}{}{#1}{#2}}

\title{Algebraic methods for computing smallest enclosing and
circumscribing cylinders of simplices}

\author{Ren\'{e} Brandenberg}
\address{Ren\'{e} Brandenberg \\
        Zentrum Mathematik\\
        Technische Universit\"at M\"unchen\\
        Boltzmannstr.~3\\
        D--85747 Gar\-ching bei M\"unchen}
\curraddr{on leave at Technische Universit\"at Wien \\
        Institut f\"ur Analysis und Technische Mathematik \\
        Wiedner Hauptstr. 8--10 \\
        A--1040 Wien}
\email{brandenb@mathematik.tu-muenchen.de}
\urladdr{http://www-m9.mathematik.tu-muenchen.de/\~{}brandenb/}

\author{Thorsten Theobald}
\address{Thorsten Theobald \\
        Zentrum Mathematik\\
        Technische Universit\"at M\"unchen\\
        Boltzmannstr.~3\\
        D--85747 Gar\-ching bei M\"unchen}
\email{theobald@mathematik.tu-muenchen.de}
\urladdr{http://www-m9.mathematik.tu-muenchen.de/\~{}theobald/}

\date{\today}
\begin{document}

\begin{abstract}
We provide an algebraic framework to compute smallest enclosing
and smallest circumscribing cylinders of simplices in Euclidean 
space $\E^n$. Explicitly, the computation of a smallest enclosing
cylinder in $\mathbb{E}^3$ is reduced to the 
computation of a smallest circumscribing cylinder.
We improve existing polynomial formulations to compute the locally
extreme circumscribing cylinders in $\E^3$ and exhibit subclasses
of simplices where the algebraic degrees can be further reduced.
Moreover,
we generalize these efficient formulations to the $n$-dimensional case and
provide bounds on the number of local extrema.
Using elementary invariant theory,
we prove structural results 
on the direction vectors of any locally extreme circumscribing cylinder
for regular simplices.
\end{abstract}

\maketitle

\markboth{\uppercase{R.~Brandenberg and T.~Theobald}}
  {\uppercase{Smallest enclosing and circumscribing cylinders}}

\section{Introduction}

Radii (of various types) belong to the most important functionals of
polytopes and general convex bodies 
in Euclidean space $\E^n$ 
\cite{brandenberg-d2002,gritzmann-klee-dcg-92, gritzmann-klee-handbook-dcg},
and they are related to applications 
in computer vision, robotics, computational 
biology, functional analysis, and statistics (see
\cite{gritzmann-klee-mathprog-93}). 
Following the notation in \cite{brandenberg-d2002}, the 
\emph{outer $j$-radius} $R_j(\mathcal{P})$ 
of a convex body $\mathcal{C} \subset \E^n$ is the 
radius of a smallest enclosing $j$-dimensional sphere in the optimal
orthogonal projection of $\mathcal{C}$ onto a $j$-dimensional linear subspace.
Studying these radii, mainly for regular simplices and regular polytopes,
is a classical topic of convex geometry
(see \cite{bonnesen-fenchel-b34,brandenberg-regular-2002,eggleston-58,gritzmann-klee-dcg-92}).

 From the computational point of view, most of the existing algorithms for 
computing these radii focus on approximation~\cite{chan-2002, 
har-peled-varadarajan-02}.
A major reason is that exact computations lead to algebraic problems of high
degree, even for computing, say, the outer $(n{-}1)$-radius in $\E^n$
(already if $n=3$). 
However, since some approaches for computing radii of general polytopes 
consider the computation of a smallest enclosing cylinder of a simplex 
as a black box within a larger computation \cite{aas-99, ssty-2000}, 
these core problems on simplices are of fundamental importance.

Recently, the authors of~\cite{dmpt-2001} demonstrated that using 
their state-of-the-art numerical polynomial solvers, various problems
related to cylinders in $\E^3$ can be solved rather efficiently.
In particular, the authors give a polynomial formulation for the 
smallest circumscribing cylinder of a simplex in $\E^3$, whose 
B\'{e}zout 
number -- the product of the degrees of the polynomial
equations -- is 60. However, these equations contain
certain undesired solutions with multiplicity~4, and as a result of
these multiplicities the computation times (using state-of-the-art
numerical techniques) are about a factor 100 larger than those of 
similar problems in which all solutions occur with multiplicity~1.

Here, we provide a general algebraic framework for computing
smallest enclosing and circumscribing cylinders of simplices in $\E^n$. First 
we reduce the computation of a smallest enclosing cylinder in $\E^3$ 
to the computation of a smallest circumscribing cylinder, thus combining
these two problems. Then we investigate smallest circumscribing cylinders
of simplices in $\E^3$. We improve the results of~\cite{dmpt-2001} 
by providing a polynomial formulation for the locally extreme cylinders, 
whose B\'{e}zout 
bound is 36 and whose solutions generically have multiplicity one.
Our formulations use techniques from the paper \cite{MPT01} which
studies the lines simultaneously tangent to four unit spheres. 
These techniques also enable us to present classes of simplices 
for which the algebraic degrees in computing the smallest circumscribing
cylinder can be considerably reduced. 

Then, in Section~\ref{se:radiindim}, we give a generalization of our 
approach to smallest circumscribing cylinders of a simplex in $\E^n$.
Based on this formulation we give bounds on the number of 
locally extreme cylinders based on the B\'{e}zout number.
Since this bound is not tight, we provide better bounds for small
dimensions; these bounds are based on mixed volume computations and
Bernstein's Theorem.
Moreover, we study in detail the locally extreme circumscribing cylinders
of a regular simplex in $\E^n$. To exploit many symmetries in the
analysis, we provide a formulation based on symmetric polynomials.
Using elementary invariant theory we show that the direction vector
of every locally extreme circumscribing cylinder has at most three 
distinct values in its components. With this result we can illustrate
our combinatorial results on the number of solutions for general simplices.

As a byproduct of our computational studies, 
we discovered a subtle but severe mistake 
in the paper~\cite{weissbach-83-erroneous} on the explicit determination of 
the outer $(n{-}1)$-radius for a regular simplex in $\E^n$, 
thus completely invalidating the proof given there.
In the appendix we give a description of that flaw, including some 
computer-algebraic calculations illustrating it.

\section{Preliminaries and background}

\subsection{$j$-radii and cylinders}
Throughout the paper we work in Euclidean space $\E^n$, 
i.e., $\mathbb{R}^n$ with the usual scalar product 
$x \cdot y = \sum_{i=1}^3 x_i y_i$ and norm $||x|| = (x \cdot x)^{1/2}$.
We write $x^2$ for $x \cdot x$.

A \emph{$j$-flat} is an affine subspace of dimension $j$. 
For a convex polytope $\mathcal{P} \subset \E^n$ (or a finite point set
$\mathcal{P} \subset \E^n$) 
and a $j$-flat $E$, we consider
\[
  \mathcal{RD}(\mathcal{P},E) := \max_{p \in \mathcal{P}} \text{dist}(p, E),
\]
where $\text{dist}(p,E)$ denotes the Euclidean distance from $p$ to $E$.
The \emph{outer $j$-radius} of $\mathcal{P}$ is
\[
  R_j(\mathcal{P}) := \min_{E \text{ is an $(n{-}j)$-flat}} \mathcal{RD}(\mathcal{P}, E) \, .
\]

The choice of the indexing in the $j$-radius stems from the 
fact that it measures
the radius of an enclosing $j$-dimensional sphere in the optimal orthogonal
projection of $\mathcal{P}$ onto a $j$-dimensional linear subspace
(cf.\ \cite{brandenberg-d2002, gritzmann-klee-dcg-92}).

One of the most natural representatives of this class is the one with
$j=2$, $n=3$, i.e., the smallest enclosing (circular)
cylinder of a polytope.
In $\E^n$, we define a cylinder to be a set of the form
\[
  \text{bd}(\ell + \rho \mathbb{B}^n),
\]
where $\ell$ is a line in $\E^n$, $\mathbb{B}^n$ denotes the unit ball,
$\rho > 0$, the addition denotes the Minkowski sum, and $\text{bd}(\cdot)$
denotes the boundary of a set. We say that $P$ \emph{can be enclosed}
in a cylinder $\mathcal{C}$ if $P$ is contained in the convex hull of
$\mathcal{C}$.
Thus the outer $(n{-}1)$-radius gives the radius of the smallest enclosing
cylinder of a polytope.

A simplex in $\E^n$ is the convex hull of $n+1$ affinely independent
points. 
An enclosing cylinder $\mathcal{C}$ of a simplex $\mathcal{P}$ is called a 
\emph{circumscribing} cylinder of $\mathcal{P}$ if all the
vertices of $\mathcal{P}$ are contained in (the hypersurface) $\mathcal{C}$.

\subsection{Smallest circumscribing cylinders and smallest enclosing cylinders}

The following statement connects the computation of a smallest enclosing
cylinder of a polytope with the computation of a smallest circumscribing 
cylinder of a simplex.\footnote{We remark that a similar statement has already been used 
in~\cite{ssty-2000}, but the manuscript referenced there does not 
contain a complete proof.}

\begin{thm}
\label{th:smallestenclosing}
Let $\mathcal{P} = \{p_1, \ldots, p_m\}$ be a set of $m \ge 4$ points
in $\E^3$, not all collinear. If $\mathcal{P}$ can be enclosed in a
circular cylinder $\mathcal{C}$ of radius $r$,
then there exists a circular cylinder $\mathcal{C'}$ of radius $r$
enclosing all elements of $\mathcal{P}$ 
such that the surface $\mathcal{C'}$ passes through
\begin{enumerate}
\item[(i)] at least four non-collinear points of $\mathcal{P}$, or
\item[(ii)] three non-collinear points of $\mathcal{P}$, and the axis 
  $\ell$ of $\mathcal{C'}$ is contained in
  \begin{enumerate}
    \item the cylinder naturally defined by spheres of 
       radius $r$ centered at two of these points;
    \item the double cone 
       naturally defined by spheres of radius $r$ centered at two of 
       these points
       (and these spheres are disjoint);
    \item or the set of lines which are tangent to the two spheres of radius 
       $r$ centered at two these points and which 
       are contained in the plane equidistant from these 
       points (and the spheres are non-disjoint).
  \end{enumerate}
\end{enumerate}
Moreover, $\mathcal{C}$ can be transformed into $\mathcal{C'}$ by a continuous
motion.
\end{thm}

\ifpictures
\begin{figure}[tb]

\[
  \begin{array}{c@{\hspace*{-1cm}}c}
    \includegraphics[scale=0.8]{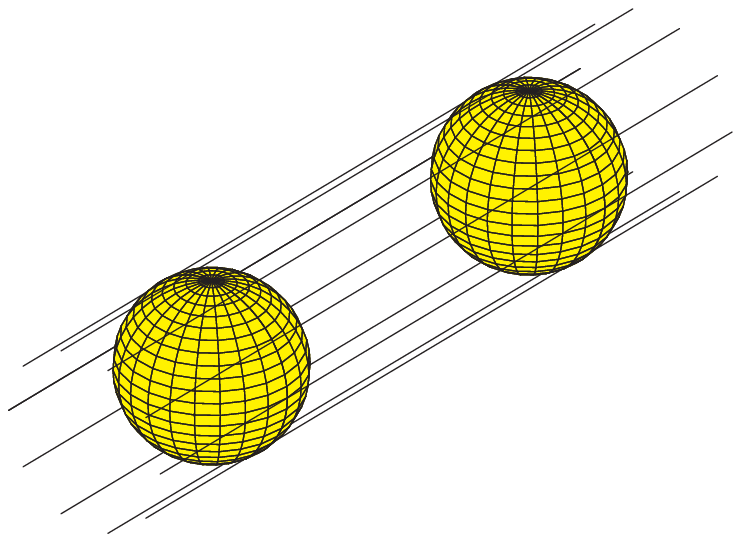} &
    \includegraphics[scale=0.8]{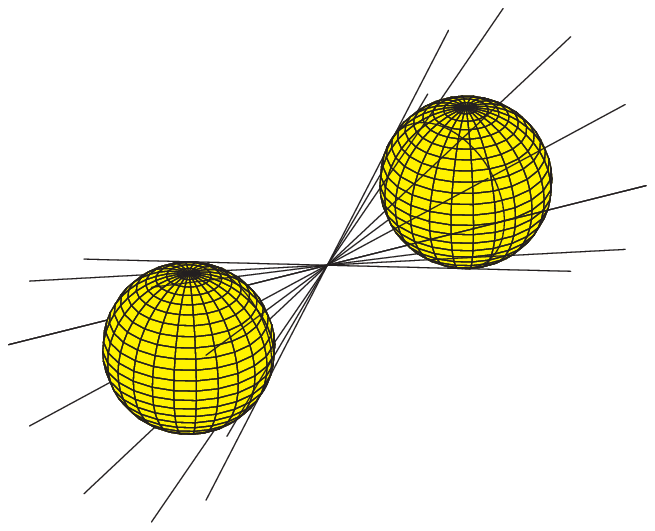} \\
    \mbox{(a) Cylinder} &
    \mbox{\hspace*{2cm}(b) Double cone with apex $(a/2,0,0)^T$}
  \end{array}
\]
\vspace*{-0.3cm}

\caption{Extreme situations of the set of hyperboloids for disjoint
spheres}
\label{fi:cylindercone}
\end{figure}
\fi

Figures~\ref{fi:cylindercone} and~\ref{fi:hyperboloid1} 
visualize the three geometric properties in the second possibility.

\ifpictures
\begin{figure}[tb]

\[
  \begin{array}{c@{\hspace*{-1cm}}c}
    \includegraphics[scale=0.8]{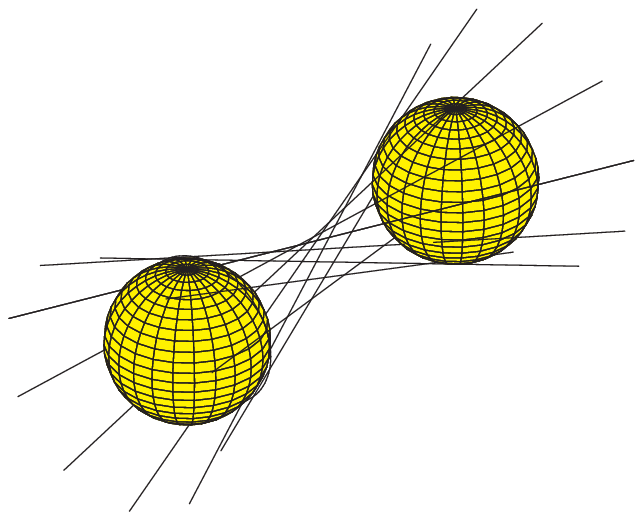} &
    \includegraphics[scale=0.8]{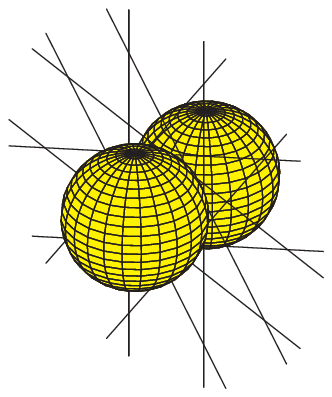} \\
    \mbox{(a) Hyperboloid for $0 < x_h < 2r^2/a$} &
    \mbox{\hspace*{2cm} (b) Degenerated hyperboloid for $x_h=a/2$}
  \end{array}
\]
\vspace*{-0.3cm}

\caption{The left figure shows a general situation for disjoint spheres;
the right figure shows an extreme situation for non-disjoint spheres}
\label{fi:hyperboloid1}
\end{figure}
\fi

Since the second possibility in Theorem~\ref{th:smallestenclosing} 
characterizes the possible special cases, this lemma in particular
reduces the computation of a smallest enclosing cylinder
of a simplex in $\E^3$ to the computation of a smallest circumscribing cylinder
of a simplex. Namely, it suffices to compute the smallest circumscribing
cylinder (corresponding to case (i)) as well as the smallest enclosing
cylinders whose axes satisfies one of the condition in (ii); the
latter case gives a constant number of problems of smaller algebraic
degree (since the positions of the axes are very restricted).

\begin{rem}
Before we start with the proof, we remark that 
Theorem \ref{th:smallestenclosing} and its
different cases show a quite similar behaviour as the well known
statement that the (unique) circumsphere of a simplex touches all its
vertices, or one of its great $(n{-}1)$-circles is the circumsphere of one
of the $(n{-}1)$-faces of the simplex (see \cite[p.~54]{bonnesen-fenchel-b34}).
\end{rem}

In the proof we will apply the following geometric equivalence. A
point $x \in \E^3$ is enclosed in a cylinder with axis $\ell$
if and only if $\ell$ is a transversal of the sphere with radius $r$ 
centered at $x$ (i.e., $\ell$ is a line intersecting the sphere).

\begin{proof}[Proof of Theorem \ref{th:smallestenclosing}] Let $\mathcal{C}$ be a cylinder with axis $\ell$ and radius
$r$ enclosing $\mathcal{P}$. Then, denoting by $S_i := S(p_i, r)$ the sphere
with radius $r$ centered at $p_i$, $\ell$ is a common transversal to
$S_1, \ldots, S_m$.
By continuously translating and rotating $\mathcal{\ell}$, we can assume 
that $\ell$ is tangent to two of the spheres, say $S_1$ and $S_2$.
Further, by changing coordinates, we can assume that $S_1$ and $S_2$ have the
form $S_1 = S((0,0,0)^T,r)$, $S_2 = S((a,0,0)^T,r)$ for some $a > 0$.

The set of lines tangent to two spheres of radius $r$ constitutes a 
set of hyperboloids (see, e.g., \cite{coxeter-b61,hcv-b32}).
Moreover, any of these hyperboloids touches the sphere $S_1$ on a circle lying 
in a hyperplane parallel to the $yz$-plane.
Hence, the set of hyperboloids can be parametrized by the $x$-coordinate
of this hyperplane which we denote by $x_h$.

If $S_1 \cap S_2 = \emptyset$ then the boundary values are $x_h=0$ 
and $x_h = 2r^2/a$. These two extreme situations yield a cylinder
and a double cone with apex $(a/2,0,0)^T$,
respectively (see Figure~\ref{fi:cylindercone}).
For $0 < x_h < 2r^2/a$ we obtain a hyperboloid of one sheet 
(see Figure~\ref{fi:hyperboloid1}(a)).

If $S_1 \cap S_2 \neq \emptyset$ then the boundary values are $x_h = 0$ and
$x_h = a/2$. Here, for $0 < x_h < a/2$ we obtain hyperboloids of one sheet,
too. For $x_h = a/2$ the hyperboloid degenerates to a set of tangents which
are tangents to the circle with radius
$r_c = \sqrt{4r^2-a^2}$ in the hyperplane $x = a/2$ 
(see Figure~\ref{fi:hyperboloid1}(b)).

Let $x_{h,0}$ be the parameter value of the hyperboloid containing the
line $\ell$.
The tangent to $S_1$ and $S_2$ is contained in the hyperboloid with some
parameter value $x_{h,0}$. By decreasing the parameter $x_h$ starting
from $x_{h,0}$ the hyperboloid changes its shape towards the cylinder 
around $S_1$ and $S_2$. Let $x_{h,1}$ be the infimum of all
 $0 \le x_h < x_{h,0}$ such that the hyperboloid does not contain a 
generating line tangent to some other sphere $S(p_i, r)$ for some 
$3 \le i \le m$.
If $x_{h,1} = 0$, then by choosing any point of $\mathcal{P}$
not collinear to
$p_1$ and $p_2$ we are in case (ii) (a).

If $x_{h,1} > 0$ then let $p_3$ be the corresponding point.
Let $T(S_1,S_2,S_3)$ denote the set of lines simultaneously tangent to
$S_1$, $S_2$, and $S_3$.
Now let $x_{h,2}$ be the infimum of all $0 \le x_h < x_{h,0}$ such that
there exists a continuous function 
$\ell : (x_{h,2},x_{h,1}) \to T(\{S_1,S_2,S_3\})$
with $\ell(x_h)$ lying on the hyperboloid with parameter $x_h$.
Since the spheres are compact, the infimum is a minimum.
If $x_{h,2} > 0$ then one of three hyperboloids involved by the
three pairs of spheres must be one of the extreme hyperboloids in 
that situation and we are in cases (ii) (a), (b), or (c).
If $x_{h,2} = 0$ then we distinguish between two possibilities. Either
during this process we also reached a tangent to some other sphere
$S(p_i,r)$ for some $4 \le i \le m$; in this case we are in case (i).
Or during the transformation all the points $p_4, \ldots, p_m$ are
enclosed in the cylinder with axis $\ell$ and radius $r$, but none
of them is contained in it. Then
we arrive at situation (ii) (a).
\end{proof}

\section{Computing the smallest circumscribing cylinders of a simplex
in $\E^3$\label{se:smallestcircumscribingr3}}

So far, we have seen how to reduce the computation of a smallest
enclosing cylinder of a simplex in $\E^3$ to the computation
of a smallest circumscribing cylinder. In order to apply algebraic 
methods to compute a smallest circumscribing cylinder, there are many
different ways to formulate that problem in terms of polynomial equations.
It is well-known that the computational costs of solving a system of
polynomial equations are mainly dominated by the B\'{e}zout number
(= product of the degrees) and the mixed volume (the latter one is discussed
in Section~\ref{se:radiindim}). 
See~\cite{clo-b96,clo-b98,sturmfels-cbms} 
for comprehensive introductions and the state-of-the-art.
Hence, it is an essential task to find the right formulations.
Moreover, we are interested in simplex classes for which the degrees
can be further reduced.

\subsection{General simplices in $\E^3$\label{su:generalsimplices}}
In the proof of~\cite[Theorem~6]{dmpt-2001}, 
a polynomial formulation is given to compute the smallest enclosing cylinder
of a simplex in $\E^3$. This formulation describes the problem by
three equations in the direction vector $v = (v_1,v_2,v_3)^T$ of the line,
one of them normalizing the direction vector $v$ by
\begin{equation}
  \label{eq:normcond}
  v_1^2 + v_2^2 + v_3^2 \ = \ 1 \, .
\end{equation}
The equations are of degree 10, 3, and 2, respectively, thus giving
a B\'{e}zout number of 60. However, as pointed out in that paper,
some of the solutions to that system are artificially introduced by
the formulation and occur with higher multiplicity, and there
are only 18 really different solutions. Even more severely,
in the experiments in that paper 
(using \textsc{Synaps}, a state-of-the-art software 
for numerical polynomial computations),
the numerical treatment of these multiple solutions needs much time,
roughly a factor 100 compared to similar systems without multiple solutions.

Here, we present an approach, which reflects the true algebraic 
bound of 18. Namely, we give a polynomial formulation with
B\'{e}zout bound 36 in which every solution generically has
multiplicity one. The additional factor 2 just results from the
fact that due to the normalization condition~\eqref{eq:normcond} 
every solution $v$ also implies that $-v$ is a solution as well.

Our framework is based on~\cite{MPT01} in which the lines simultaneously
tangent to four unit spheres are studied.
A line in $\E^3$ is represented by a point $u \in \E^3$
lying on the line and a direction vector 
$\bv\in \E^{3}$ with $v^2 = 1$.
We can make $u$ unique by requiring that $u \cdot v = 0$.
A line $\ell = (u,\bv)$ has Euclidean distance $r$ from a point 
$p \in \E^3$ if and only if the quadratic equation 
$(u +tv - p)^2 = r^2$ has a solution of multiplicity two.
This gives the condition
$$
   \frac{\left(\bv\cdot(u-p)\right)^2}{\bv^2} - (u-p)^2\ + r^2\ =\ 0 \,.
$$
Expanding this equation yields
\begin{equation}\label{eq:tanSphere}
  v^2 u^2- 2 v^2 u \cdot p + v^2 p^2 - (v \cdot p)^2-r^2\bv^2\ =\ 0\,. \end{equation}
Rather than using $v^2 = 1$ to further simplify this equation, we prefer
to keep the homogenous form, in which all terms are of degree~4.

Now let $p_1, \ldots, p_4$ be the affinely independent vertices of the
given simplex.
Without loss of generality we can choose
$p_4$ to be the located in the origin.
Then the remaining points span $\E^3$.
Subtracting the equation for the point in the origin from the equations for
$p_1$, $p_2$, $p_3$ gives the following program to compute the 
square of the radius of the minimal circumscribing cylinder.
 \begin{equation}\label{eq:origintransform}
  \begin{array}{rrcl}
   & \min u^2 \\
   \text{s.t.} & u \cdot\bv & = & 0 \, ,\\
      & 2 \bv^2 u \cdot p_i & = & \bv^2 p_i^2 - (v \cdot p_i)^2 \, ,
      \qquad 1 \le i \le 3\, , \\
      & v^2 & = & 1 \, .
  \end{array}
 \end{equation}

We remark that the set of admissible solutions is nonempty; a proof 
of that statement (for general dimension) is contained in 
Section~\ref{se:radiindim}.

Since the points $p_1, p_2, p_3$ are linearly independent,
the matrix $M := (p_1, p_2, p_3)^T$ is invertible, and we can
solve the equations in the bottom line of~\eqref{eq:origintransform} for~$u$:
 \begin{equation}\label{eq:pequation}
  u \ =\ \frac{1}{2\bv^2} M^{-1}
    \left( \begin{array}{c}
      \bv^2 p_1^2 - (v \cdot p_1)^2\\
      \bv^2 p_2^2 - (v \cdot p_2)^2\\
      \bv^2 p_3^2 - (v \cdot p_3)^2
    \end{array} \right).
 \end{equation}
Now substitute this expression for $u$ into the the objective function and
into the first constraint
of the system~(\ref{eq:origintransform}). After setting $v^2 = 1$
in the denominator
of the first constraint, this gives a homogeneous cubic equation which we 
denote by $g_1(v_1,v_2,v_3) = 0$. Hence, we arrive at the
following polynomial optimization formulation in terms of the
variables $v_1$, $v_2$, and $v_3$. 
\begin{equation}
  \label{eq:mincyl3d}
  \begin{array}{rrcl}
  & \multicolumn{3}{l}{\min\left( \frac{1}{2} M^{-1}
    \left( \begin{array}{c}
      \bv^2 p_1^2 - (v \cdot p_1)^2 \\
      \bv^2 p_2^2 - (v \cdot p_2)^2 \\
      \bv^2 p_3^2 - (v \cdot p_3)^2
    \end{array} \right)
    \right)^2} \\ [2ex]
  \text{s.t.} & g_1(v_1,v_2,v_3) & = & 0 \, , \\
  &  g_2(v_1,v_2,v_3) := v^2 - 1 & = & 0 \, .
  \end{array} 
\end{equation}
Note that the objective function is a homogeneous 
polynomial of degree 4. We denote this polynomial by $f$.

Using Lagrange multipliers $\lambda_1$ and $\lambda_2$, a necessary
local optimality condition is
\begin{equation}
\label{eq:lagrange1}
  \text{grad } f= \lambda_1 \text{grad } g_1 + \lambda_2 \text{grad } 
g_2 \, .
\end{equation}
By thinking of an additional factor $\lambda_0$ before $\text{grad }f$ 
and considering~\eqref{eq:lagrange1} as a system of linear equations in
$\lambda_0$, $\lambda_1$, $\lambda_2$, we see that 
if~\eqref{eq:lagrange1} is satisfied for some vector $v$
then the determinant
\begin{equation}
\label{eq:optimalitydet}
  \det \left( \begin{array}{ccc}
    -\frac{\partial f}{\partial v_1} & \frac{\partial g_1}{\partial v_1} & 
     \frac{\partial g_2}{\partial v_1} \\
    -\frac{\partial f}{\partial v_2} & \frac{\partial g_1}{\partial v_2} & 
     \frac{\partial g_2}{\partial v_2} \\
    -\frac{\partial f}{\partial v_3} & \frac{\partial g_1}{\partial v_3} &
     \frac{\partial g_2}{\partial v_3}
  \end{array} \right) 
\end{equation}
vanishes.

\begin{lemma}
(a) Any direction vector $(v_1, v_2, v_3)^T \in \E^3$ 
of the axis of a locally 
extreme circumscribing cylinder satisfies the polynomial
system~\eqref{eq:mincyl3d}.
If this system has only finitely many solutions then this number is
bounded by~36.

(b) For a generic simplex the number of solutions is indeed finite, 
and all solutions have multiplicity one.
\end{lemma}

\begin{proof}
Let $v$ be the direction vector of an axis of a locally
extreme circumscribing cylinder. 
Then $v$ satisfies the first constraint of~\eqref{eq:mincyl3d}, and the
determinant~\eqref{eq:optimalitydet} vanishes. Since these are homogeneous
equations of degree~3 and 6, respectively, B\'{e}zout's Theorem implies
that in connection with $v^2=1$ we obtain at most 36 isolated solutions.

For the second statement it suffices to check that for one specific
simplex there are only finitely many solutions and that all solutions are 
pairwise distinct.
\end{proof}

\subsection{Special simplex classes in $\E^3$}

In this section, we investigate conditions under which the degree 
of the resulting
equations is reduced. Moreover, we show that for the equifacial simplex,
the minimal circumscribing radius can be computed quite easily.

We use the following classification from \cite{MPT01,schaal-85}.

\begin{prop}
\label{pr:charactreducible}
Let $T$ be a simplex in $\E^3$ with vertices $p_1, \ldots, p_4$.
The polynomial $g_1$ in the cubic equation factors into a linear polynomial 
and an irreducible quadratic polynomial 
if and only if the four faces of $T$ can
be partitioned into two pairs of faces $\{F_1,F_2\}$, $\{F_3, F_4\}$
with $\textrm{area}(F_1) = \textrm{area}(F_2) \neq
\textrm{area}(F_3) = \textrm{area}(F_4)$.
Moreover, $g_1$ factors into three linear terms if and only if
the areas of all four faces of $T$ are equal.
\end{prop}

First let us consider the case where $g_1$ decomposes 
into a linear polynomial and an irreducible quadratic polynomial. 
By optimizing separately over the linear and the quadratic 
constraint, the degrees of our equations are smaller than for the general
case. Namely, analogously to the derivation in 
Section~\ref{su:generalsimplices},
for the quadratic constraint we obtain a B\'{e}zout bound of 
\[
  (3+1+1) \cdot 2 \cdot 2 = 20 \, ,
\]
and for the linear constraint we obtain
\[
  (3+0+1) \cdot 1 \cdot 2 = 8 \, .
\]
Thus, we can conclude:

\begin{lemma}
If the four faces of the simplex can be partitioned into two pairs of
faces $\{F_1,F_2\}$, $\{F_3,F_4\}$ 
with $\textrm{area}(F_1) = \textrm{area}(F_2) \neq
\textrm{area}(F_3) = \textrm{area}(F_4)$
then there are at most 28 isolated local
extrema for the mimimal circumscribing cylinder. They can be computed
from two polynomial systems with B\'{e}zout numbers 20 and 8, respectively.
\end{lemma}

\bigskip

\noindent
{\bf Equifacial simplices.} 
A simplex in $\E^3$ is called \emph{equifacial} if all four faces have the 
same area. By Proposition~\ref{pr:charactreducible}, for an equifacial simplex
the cubic polynomial $g_1$
factors into three linear terms. Hence, 
we obtain at most $3 \cdot 8 = 24$ local extrema.
Somewhat surprisingly, using a characterization 
from~\cite{theobald-dimacs-2001}, 
it is even possible to compute smallest circumscribing cylinder 
of an equifacial simplex esentially without any algebraic computation.

Namely, it is well-known that the vertices of an equifacial
simplex $T$ can be regarded as four pairwise non-adjacent
vertices of a rectangular box
(see, e.g., \cite{kupitz-martini-99}).
Hence, there exists a representation
$p_1=(w_1, w_2,w_3)^T$,
$p_2=(w_1,-w_2,-w_3)^T$,
$p_3=(-w_1,w_2,-w_3)^T$,
$p_4=(-w_1,-w_2,w_3)^T$ with
$w_1, w_2, w_3 > 0$. 

Assuming without loss of generality $v^2 = 1$, 
(\ref{eq:tanSphere}) gives
\begin{equation}
  \label{eq:equifacialstart}
   (v \cdot p_i)^2 + 2 u \cdot p_i = \sum_{j=1}^3 w_j^2
  + u^2 - r^2, \qquad 1 \le i \le 4 \, .
\end{equation}
Subtracting these equations pairwise gives
\[
  4 (w_2 u_2 + w_3 u_3) = - 4 (w_1 w_3 v_1 v_3 +
  w_1 w_2 v_1 v_2)
\]
(for indices 1, 2) and analogous equations, so that
\[
  w_1 u_1 = - w_2 w_3 v_2 v_3, \quad
  w_2 u_2 = - w_1 w_3 v_1 v_3, \quad
  w_3 u_3 = - w_1 w_2 v_1 v_2.
\]
Since $u \cdot v = 0$, this yields $v_1 v_2 v_3 = 0$.
Without loss of generality we can assume $v_1 = 0$. In this case,
\[
  u = \left( - \frac{w_2 w_3}{w_1} v_2 v_3, 0, 0 \right)^T.
\]
So we can express~(\ref{eq:equifacialstart}) in terms of the
direction vector $v$,
\[
  w_2^2 v_2^2 + w_3^2 v_3^2 = \sum_{j=1}^3 w_j^2 +
  \left( - \frac{w_2 w_3}{w_1} v_2 v_3 \right)^2
  - r^2,
\]
which, by using $v_2^2 + v_3^2 = 1$, gives
\begin{equation}
\label{eq:r2intermsofs2}
  r^2 \ = \ 
  - \frac{w_2^2 w_3^2}{w_1^2} v_2^4 -
  \left( w_2^2 - w_3^2 -
   \frac{w_2^2 w_3^2}{w_1^2} \right) v_2^2
  + w_1^2 + w_2^2 \, .
\end{equation}
Thus, by computing the derivative of this expression $r^2 = r^2(v_2)$
and taking into account
the three cases $v_i = 0$, we can reduce the
computation of the minimal circumscribing cylinder to solving
three univariate equations of degree~3.
However, we can still do better.
Substitute $z_2 := v_2^2$, and let $\rho$ be the expression for $r^2$ in terms
of $z_2$,
\[
   \rho(z_2) \ = \ 
  - \frac{w_2^2 w_3^2}{w_1^2} z_2^2 -
  \left( w_2^2 - w_3^2 -
   \frac{w_2^2 w_3^2}{w_1^2} \right) z_2
  + w_1^2 + w_2^2 \, .
\]
Since the second derivative of that quadratic function is negative, 
$\rho(z_2)$ is a concave function. Hence, within the interval
$z_2 \in [0,1]$, the minimum is attained at one of the boundary
values $z_2 \in \{0,1\}$. Consequently, two of the components of 
$(v_1,v_2,v_3)^T$ must be zero and therefore $v$ is perpendicular to
two opposite edges. Since the latter geometric characterization is 
independent of our specific choice of coordinates, we can conclude:

\begin{lemma}
\label{eq:equifacialminradius}
If all four faces of the simplex $T$ have the same area then the axis
of a minimum circumscribing cylinder is perpendicular to 
two opposite edges.
\end{lemma}

Hence, for an equifacial simplex it suffices to investigate
the cross products of the three pairs of opposite edges (equipped
with an orientation), 
and we do not need to solve a system of polynomial equations at all.

In order to illustrate how these three solutions relate to the 
18 solutions of the general approach above,
we consider the regular simplex in $\E^3$.
In the general approach, as already pointed out in \cite{dmpt-2001}, 
the six edge directions $p_i p_j$ ($1 \le i < j \le 4$) all have
multiplicity~1, and each of the three directions in
Lemma~\ref{eq:equifacialminradius},
$p_1 p_2 \times p_3 p_4$, 
$p_1 p_3 \times p_2 p_4$, 
$p_1 p_4 \times p_2 p_3$, have multiplicity 4.

\section{Smallest circumscribing cylinders in higher dimensions\label{se:radiindim}}

In Section~\ref{se:smallestcircumscribingr3} we have given polynomial 
formulations with small B\'{e}zout number for computing smallest 
circumscribing cylinders of a simplex in $\E^3$. Using the 
characterization in~\cite{sottile-theobald-ndim} of lines
simultaneously tangent to $2n{-}2$ spheres in $\E^n$, we generalize these
formulations to smallest circumscribing cylinders of a simplex 
in $\E^n$, $n \ge 2$.
Analogous to the three-dimensional case let $p_1, \ldots, p_{n+1}$ be
the affinely independent 
vertices of the simplex in $\E^n$, and let
$p_{n+1}$ be located in the origin.

First note that~\eqref{eq:origintransform} also holds in general dimension $n$
if we replace the index $3$ by the index $n$.
Since the points $p_1, \ldots, p_n$ are linearly independent,
the matrix $M := (p_1, \ldots, p_n)^T$ is invertible, and we can solve 
for $u$:
\begin{equation}\label{eq:pequationdimn}
  u \ =\ \frac{1}{2\bv^2} M^{-1}
    \left( \begin{array}{c}
      \bv^2 p_1^2 - (v \cdot p_1)^2 \\
      \vdots \\
      \bv^2 p_n^2 - (v \cdot p_n)^2
    \end{array} \right).
 \end{equation}

Hence, by generalizing the formulation for the three-dimensional case, 
we obtain the program
\begin{equation}
  \label{eq:optformulationndim}
  \begin{array}{rrcl}
  & \multicolumn{3}{r}{\min\left( \frac{1}{2} M^{-1}
    \left( \begin{array}{c}
      \bv^2 p_1^2 - (\bv \cdot p_1)^2 \\
      \vdots \\
      \bv^2 p_n^2 - (\bv \cdot p_n)^2
    \end{array} \right)
    \right)^2} \\ [2ex]
  \text{s.t.} & g_1(v_1, \ldots, v_n) & = & 0 \, , \\
  & g_2(v_1, \ldots, v_n) := v^2 -1 & = & 0 \, ,
  \end{array}
\end{equation}
where $g_1$ denotes the cubic equation as before.
In order to show that set of admissible solutions for our optimization
problem is nonempty, we record the following result.

\begin{lemma}
For any simplex in $\E^n$ the $\binom{n+1}{2}$ edge directions of the
simplex are direction vectors of circumscribing cylinders.
\end{lemma}

\begin{proof}
Since the edge directions $p_i - p_j$ have a simple description in
the basis $p_1, \ldots, p_n$, we express 
the cubic equation $g_1(v)=0$ in that basis.
Let $v$ be an arbitrary direction vector, and let the representation of $v$ in 
the basis $p_1, \ldots, p_n$ be
\[
  v\ =\ \sum_{i=1}^n t_i p_i \, .
\]
Further, let $p_1', \ldots, p_n'$ be a dual basis to $p_1, \ldots, p_n$;
i.e., let $p_1', \ldots, p_n'$ be defined by 
$p_i' \cdot p_j = \delta_{ij}$, where $\delta_{ij}$ denotes Kronecker's
delta function. 
By elementary linear algebra, we have $t_i = p_i' \cdot v$.

When expressing $u$ in this dual basis,
$u = \sum u_i' p_i'$, the second constraint of~(\ref{eq:origintransform}) 
gives
\[
  u_i'\ =\ \frac{1}{2 v^2}
    \left(v^2 p_i^2 - (v \cdot p_i)^2 \right) \, .
\]
Substituting this representation of $u$ into the equation $g_1(v) = 0$
gives
\[
  0\ =\ g_1(v) =\ v^2 (u \cdot v)
   \ =\ v^2 \left(\sum_{i=1}^n u_i' p_i'\right) \cdot v 
   \ =\ v^2 \sum_{i=1}^n u_i' t_i \, ,
\]
where the last step uses the duality of the bases.
Hence, we obtain the cubic equation
\[
  \frac{1}{2} \sum_{i=1}^n (v^2 p_i^2 - (v \cdot p_i)^2) t_i \ =\ 0 \, .
\]
Expressing $v$ in terms of the $t$-variables yields
\[
  \frac{1}{2} \sum_{1 \le i \neq j \le n}
    \alpha_{ij} t_i^2 t_j + 
    \sum_{1 \le i < j < k \le n} \beta_{ijk} t_i t_j t_k\ =\ 0 \, ,
\]
where
\begin{eqnarray*}
  \alpha_{ij} & = & (\text{vol}_2(p_i,p_j))^2\ =\ 
    \det \left( \begin{array}{cc}
      p_i\cdot p_i & p_i \cdot p_j  \\
       p_j\cdot p_i & p_j \cdot p_j 
  \end{array} \right), \\
  \beta_{ijk} & = &
    \det \left( \begin{array}{cc}
       p_i\cdot p_j  &  p_i \cdot p_k  \\
       p_k\cdot p_j  &  p_k \cdot p_k 
  \end{array} \right) +
    \det \left( \begin{array}{cc}
       p_i\cdot p_k  &  p_i \cdot p_j  \\
       p_j\cdot p_k  &  p_j \cdot p_j 
  \end{array} \right) \\
  & & +
    \det \left( \begin{array}{cc}
       p_j\cdot p_k  &  p_j \cdot p_i  \\
       p_i\cdot p_k  &  p_i \cdot p_i 
  \end{array} \right),
\end{eqnarray*}
and $\text{vol}_2(p_i,p_j)$ denotes the oriented area of the parallelogram
spanned by $p_i$ and $p_j$.
In terms of the $t$-coordinates, the $\binom{n+1}{2}$ edges of the simplex
are $t = e_i$, $1 \le i \le n$, and $t = e_i - e_j$, $1 \le i < j \le n$,
where $e_i$ denotes the $i$-th standard unit vector.
For all these edges, the cubic equation is satisfied.
\end{proof}

Considering Lagrange multipliers $\lambda_1$ and $\lambda_2$ yields
the following necessary optimality condition.
\begin{eqnarray}
  \text{grad } f & = & \lambda_1 \text{grad } g_1 + \lambda_2 \text{grad } 
g_2 \, , \nonumber \\
  g_1(v_1, \ldots, v_n) & = & 0 \, , \label{eq:lagrangendim} \\
  g_2(v_1, \ldots, v_n) & = & 0 \, . \nonumber
\end{eqnarray}
Since the B\'{e}zout bound of this system is 
$3^n \cdot 3 \cdot 2 = 2 \cdot 3^{n+1}$, we have:

\begin{lemma}
\label{le:boundndim}
For $n \ge 2$, the number of isolated local extrema for the minimal
circumscribing cylinder is bounded by $2 \cdot 3^{n+1}$.
\end{lemma}

This bound is not tight. Trying to reduce this upper bound of isolated 
solutions like in the three-dimensional case, we can eliminate the
linear occurrences of the Lagrange variables $\lambda_1$ and $\lambda_2$.
Generalizing~\eqref{eq:optimalitydet}, we have to consider the
vanishing of all $3 \times 3$-subdeterminants of the matrix
\begin{equation}
\label{eq:optimalitydetndim}
   \left( \begin{array}{ccc}
    -\frac{\partial f}{\partial v_1} & \frac{\partial g_1}{\partial v_1} & 
     \frac{\partial g_2}{\partial v_1} \\
    -\frac{\partial f}{\partial v_2} & \frac{\partial g_1}{\partial v_2} & 
     \frac{\partial g_2}{\partial v_2} \\
     \vdots & \vdots & \vdots \\
    -\frac{\partial f}{\partial v_n} & \frac{\partial g_1}{\partial v_n} &
     \frac{\partial g_2}{\partial v_n}
  \end{array} \right) \, . 
\end{equation}
Thus, for $n \ge 4$ we arrive at a non-complete intersection of equations 
where we have more equations than variables. Hence, we cannot apply our
B\'{e}zout bound on these systems.

However, for small dimensions we can improve Lemma~\ref{le:boundndim}
by directly working on the formulation~\eqref{eq:lagrangendim}.
In order to provide better bounds, we use well-known characterizations
of the number of zeroes of a polynomial equation by the mixed volume
of a Minkowski sum of polytopes (for an easily accessible introduction
into this topic we refer to~\cite{clo-b98}).
Here, let $\mathbb{C}^* := \mathbb{C} \setminus \{0\}$.

\begin{lemma}
\label{le:stirling1}
For $2 \le n \le 7$, the number of solutions of the 
system~\eqref{eq:lagrangendim}
in $(v_1, \ldots, v_n, \lambda_1, \lambda_2) $ $\in (\mathbb{C}^{*})^{n+2}$ is bounded
by 
\[
  6 \stirlingnumbersecond{n+1}{3} \, ,
\]
where $\stirlingnumbersecond{n}{k}$ denotes the Stirling number of the second 
kind~(see, e.g., \cite{okp-89, stanley-b97}).
\end{lemma}

The sequence $6 \stirlingnumbersecond{n+1}{3}$ starts as follows
\[
  \begin{array}{|c||c|c|c|c|c|c|c|} \hline
     n & 2 & 3 & 4 & 5 & 6 & 7 \\ \hline
     6 \stirlingnumbersecond{n+1}{3} &
         6 & 36 & 150 & 540 & 1806 & 5796 \\ \hline
  \end{array}
\]

\begin{proof}
For a polynomial 
$h = \sum_{\alpha \in \mathbb{N}_0^n} c_{\alpha} x^{\alpha} 
\in \mathbb{C}[x_1, \ldots, x_n]$, let
\[
  \mathrm{NP}(h) := \textrm{conv} \{ \alpha \in \mathbb{N}_0^n \, : \, c_\alpha \neq 0 \}
\]
denote the Newton polytope of $h$ (see, e.g., \cite[\S 7.1]{clo-b98}).
Let $h_1, \ldots, h_n$ be the polynomials of the gradient equation 
in~\eqref{eq:lagrangendim}. Further let $P_1, \ldots, P_n, Q_1, Q_2$ 
be the Newton polytopes of $h_1, \ldots, h_n, g_1, g_2$ for generic instances
of these equations.

Recall that the mixed volume 
$\textrm{MV}(P_1, \ldots, P_n, Q_1, Q_2)$ is the coefficient of the monomial 
$\lambda_1 \cdot \lambda_2 \cdots \lambda_n \cdot \mu_1 \cdot \mu_2$
in the $(n+2)$-dimensional volume 
$\textrm{Vol}_{n+2}(\lambda_1 P_1 + \ldots + \lambda_n P_n + \mu_1 Q_1 + \mu_2 Q_2)$ (which is a polynomial expression 
in $\lambda_1, \ldots, \lambda_n, \mu_1, \mu_2$).
By Bernstein's Theorem, the number of isolated common zeroes in 
$(\mathbb{C}^*)^{n+2}$ of the set of polynomials $h_1, \ldots, h_n, g_1, g_2$
is bounded aboved by 
\[
  \textrm{MV}(P_1, \ldots, P_n, Q_1, Q_2)
\]
(see \cite[Chapter~8, Theorem 5.4]{clo-b98}).
For every given $n$ this volume can be computed using software for computing
mixed volumes (see, e.g, \cite{emiris-canny-95,verschelde-phc-99}).
\end{proof}

We conjecture that for any $n \ge 2$, the number of isolated solutions 
in $(\mathbb{C}^{*})^{n+2}$ is bounded by 
$6 \stirlingnumbersecond{n+1}{3}$.

\subsection{The regular simplex in $\E^n$}

Here, we analyze the local extrema of circumscribing cylinders
for the regular simplex. Our aim is both to illustrate the
algebraic formulations given before and to relate our investigations
to classical investigations on the regular simplex in convex geometry.
In order to achieve many symmetries in the algebraic formulation, we 
use a slightly modified coordinate system that is particularly suited
for the regular simplex; these coordinates have also 
been used in~\cite{weissbach-83-correct, brandenberg-regular-2002}.

The equation $x_1 + \ldots + x_{n+1} = 1$ defines an $n$-dimensional
affine subspace in $\E^{n+1}$. Now let the 
regular simplex in this $n$-dimensional subspace be
given by the $n+1$ vertices $p_i = e_i$, where $e_i$ denotes the
$i$-th standard unit vector, $1 \le i \le n+1$. 
We consider the tangency equation~\eqref{eq:tanSphere} for the point 
$p_{n+1}$,
\[
  v^2 u^2 - 2v^2 u_{n+1} + v^2 - v_{n+1}^2 - r^2 v^2 \ = \ 0 \, .
\]
Subtracting this equation from the equation for $p_i$, $1 \le i \le n$,
yields
\[
  2 v^2 (u_i - u_{n+1}) \ = \ - (v_i^2 - v_{n+1}^2) \, , \quad 1 \le i \le n \, .
\]
Moreover, the embedding into the hyperplane $\sum_{i=1}^{n+1} x_i = 1$ 
implies $\sum_{i=1}^{n+1} u_i = 1$. In order to solve these $n+1$ equations
for $u$, let $M$ be the $(n+1) \times (n+1)$-matrix whose $i$-th 
row contains the
vector $e_i^T - e_{n+1}^T$ and whose $n$-th row is $(1,1, \ldots, 1)$.
Since $M$ is invertible, we obtain
\begin{equation}
  \label{eq:uregularsimplex}
  u = \frac{1}{2 v^2} M^{-1} \left( \begin{array}{c}
    - (v_1^2 - v_{n+1}^2) \\
    \vdots \\
    - (v_n^2 - v_{n+1}^2) \\
    2 v^2
    \end{array} \right) \, .
\end{equation}
As before, substituting this expression into $u \cdot v = 0$ and setting
$v^2=1$ in the denominator gives a cubic equation $g_1(v) = 0$. Hence,
we obtain the following optimization problem. Here, the objective function
$f$ stems from the condition for the vertex $p_{n+1}$, and the condition
$\sum_{i=1}^{n+1} v_i = 0$ comes from the embedding.
\begin{equation}
\label{eq:optregular}
\begin{array}{rrcl}
& \multicolumn{3}{l}{\min u^2 - 2 u_{n+1} + 1 - v_{n+1}^2} \\
  \text{s.t.} & g_1(v_1, \ldots, v_{n+1}) & = & 0 \, , \\
  & \sum\limits_{i=1}^{n+1} v_i & = & 0 \, , \\
  & v^2 & = & 1 \, .
  \end{array}
\end{equation}

First we record that the functions $f$ and $g_1$ are symmetric polynomials
in the variables $v_1, \ldots, v_{n+1}$. In order to show this,
let $\sigma_1, \ldots, \sigma_{n+1}$ be the elementary symmetric
functions in $v_1, \ldots, v_{n+1}$,
\begin{eqnarray*}
  \sigma_1 & = & v_1 + \ldots + v_{n+1} \, , \\
  & \vdots & \\
  \sigma_k & = & \sum_{1 \le i_1 < \ldots < i_k \le n+1}
    v_{i_1} v_{i_2} \cdots v_{i_k} \, , \\
  & \vdots & \\
  \sigma_{n+1} & = & v_1 v_2 \cdots v_{n+1} \\
\end{eqnarray*}
(see, e.g., \cite{clo-b96,sturmfels-b93}). By providing explicit
expressions for $f$ and $g_1$ as polynomials
in the elementary symmetric polynomials $\sigma_1, \ldots, \sigma_{n+1}$,
the symmetry of $f$ and $g_1$ follows. More precisely, we obtain: 

\begin{lemma}
\label{le:fg1symmetric}
The quartic polynomial $f(v_1, \ldots, v_{n+1})$ and the cubic 
polynomial $g_1(v_1, \ldots, $ $v_{n+1})$ are symmetric
polynomials in the variables $v_1, \ldots, v_{n+1}$. In terms of the
elementary symmetric functions, $f$ results in
\[
  f = \frac{1}{4(n+1)} \left(
      n \sigma_1^4 - 4 n \sigma_1^2 \sigma_2 
      + 2 (n-1) \sigma_2^2 - 4 \sigma_1^2 + 8 \sigma_2 
      + 4n \right) + \sigma_1 \sigma_3 - \sigma_4 \, ,
\]
and the homogeneous polynomial $g_1$ results in
\[
 g_1 = \frac{1}{2(n+1)} 
   \left(-(n-2) \sigma_1^3 + 3(n-1) \sigma_1 \sigma_2 \right) - \frac{3}{2} \sigma_3
   \, .
\]
\end{lemma} 

Since $\sigma_1 = 0$ and $\sum_{i=1}^{n+1} {v_i}^2 = \sigma_1^2 - 2 \sigma_2$,
we can also deduce the following formulation of our optimization problem:

\begin{cor}
\label{co:invarcor}
Finding the critical values of the minimization problem~\eqref{eq:optregular} 
is equivalent to finding the critical values $(v_1, \ldots, v_{n+1})^T$ of
the maximization problem
\begin{eqnarray}
\max~ \sigma_4 \nonumber \\
  \text{\rm s.t.} \quad \sigma_1 & = & 0 \, , \nonumber \\
  \sigma_2 & = & -\frac{1}{2} \, , \label{eq:invarcor} \\
  \sigma_3 & = & 0 \, , \nonumber
\end{eqnarray}
where $\sigma_i$ are the elementary symmetric functions in $v_1, \ldots, 
v_{n+1}$.
\end{cor}

\begin{thm}
The direction vector $(v_1, \ldots, v_{n+1})^T$ of any locally
extreme circumscribing cylinder satisfies 
$|\{ v_1, \ldots, v_{n+1} \}| \le 3$, i.e., for each solution vector
the components take at most three distinct values.
\end{thm}

\begin{proof} For $n \le 2$, the statement is trivial, so we can assume
$n \ge 3$.
Let $v$ be the direction vector of a locally extreme circumscribing
cylinder with $v^2 = 1$. Using Corollary~\ref{co:invarcor},
let $f(v) := -\sigma_4(v)$,
$g_1(v) := \sigma_3(v)$,
$g_2(v) := \sigma_2(v) - 1/2$,
and $g_3(v) := \sigma_1(v)$.
As a necessary condition for a local extremum,
for any pairwise different indices $a,b,c,d \in \{1, \ldots, n+1\}$ the 
determinant
\begin{equation}
\label{eq:3varsdet}
  \det \left( \begin{array}{cccc}
    -\frac{\partial f}{\partial v_a} & \frac{\partial g_1}{\partial v_a} &
   \frac{\partial g_2}{\partial v_a} & \frac{\partial g_3}{\partial v_a} \\ [0.5ex]
    -\frac{\partial f}{\partial v_b} & \frac{\partial g_1}{\partial v_b} &
   \frac{\partial g_2}{\partial v_b} & \frac{\partial g_3}{\partial v_b} \\ [0.5ex]
    -\frac{\partial f}{\partial v_c} & \frac{\partial g_1}{\partial v_c} &
   \frac{\partial g_2}{\partial v_c} & \frac{\partial g_3}{\partial v_c} \\ [0.5ex]
    -\frac{\partial f}{\partial v_d} & \frac{\partial g_1}{\partial v_d} &
   \frac{\partial g_2}{\partial v_d} & \frac{\partial g_3}{\partial v_d}
  \end{array} \right) 
\end{equation}
vanishes.
Since $f$, $g_1$, $g_2$, and $g_3$ are symmetric functions in the
variables $v_1, \ldots, v_{n+1}$, we can assume without loss of 
generality $a=1$, $b=2$, $c=3$, and $d=4$.
Setting $\alpha_n := \sum_{i=5}^{n+1} v_i$ and 
$\beta_n = \sum_{i=5}^{n+1} v_i^2$, we can write
\begin{eqnarray*}
  \frac{\partial g_3}{\partial v_i} & = & 1 \, , \\
  \frac{\partial g_2}{\partial v_i} & = & \sum_{\atopfrac{j=1}{j \neq i}}^4 v_j + \alpha_n \, ,  \\
  \frac{\partial g_1}{\partial v_i} & = & \sum_{\atopfrac{1 \le j < k \le 4}{j, k \neq i}} v_j v_k + \alpha_n \sum_{\atopfrac{j=1}{j \neq i}}^4 v_j + \frac{1}{2} \left( \alpha_n^2 - \beta_n \right) \,  
\end{eqnarray*}
($1 \le i \le 4$).
Moreover, since $\sigma_3(v) = 0$, we can consider 
$\sigma_3 + \frac{\partial f}{\partial v_i}$
instead of $\frac{\partial f}{\partial v_i}$. This allows to express
the resulting expression easily in terms of $\alpha_n$ and $\beta_n$.
More precisely, we obtain
\[
  \sigma_3 + \frac{\partial f}{\partial v_i} = v_i \left( \sum_{\atopfrac{1 \le j < k \le 4}{j, k \neq i}} v_j v_k + \alpha_n \sum_{\atopfrac{j=1}{j \neq i}}^4 v_j
  + \frac{1}{2} (\alpha_n^2 - \beta_n) \right) \, .
\]
Thus we can consider the determinant~\eqref{eq:3varsdet} as a polynomial
in $v_1, v_2, v_3, v_4, \alpha_n, \beta_n$.
Evaluat\-ing this $4 \times 4$-determinant $\Delta$ shows that it is independent
of $\alpha_n$, $\beta_n$ and that it factors as
\[
  \Delta = (v_1-v_2) (v_1-v_3) (v_1-v_4) (v_2-v_3) (v_2-v_4) (v_3-v_4) \, .
\]
Hence, $|\{v_1,v_2,v_3,v_4\}| \le 3$, and this holds true for any quadruple
$(a,b,c,d)$ of indices.
\end{proof}

Using this result, we illustrate the occurrence of the
Stirling numbers in Lemma~\ref{le:stirling1} for the case of a 
regular simplex. 
There are $\stirlingnumbersecond{n+1}{3}$ ways to partition the
set $V:=\{v_1, \ldots, v_{n+1}\}$ into three nonempty subsets 
$V_1$, $V_2$, $V_3$. We assume that $v_i \in V_i$, $1 \le i \le 3$, and 
that all variables within the same set take the same value.
Setting $k:=|V_1|$ and $l := |V_2|$, the formulation in 
Corollary~\ref{co:invarcor} yields the system of equations
\begin{eqnarray}
  k v_1 + l v_2 + (n+1-k-l) v_3 & = & 0 \, , \nonumber \\
  k v_1^2 + l v_2^2 + (n+1-k-l) v_3^2 & = & 1 \, , \label{eq:threevars} \\
  \sum_{\atopfrac{0 \le i_1 < i_2 < i_3 \le 3}{i_1 + i_2 + i_3 = 3}}
  \binom{k}{i_1} \binom{l}{i_2} \binom{n+1-k-l}{i_3} 
    v_1^{i_1} v_2^{i_2} v_3^{i_3}& = & 0 \, . \nonumber
\end{eqnarray}

If one of the indices $k$, $l$, or $n+1-k-l$ is zero then this system 
consists of three equations in two variables, so we do not expect any 
solutions. 
For every choice of $k$, $l$ corresponding to a partition into
nonempty subsets, we obtain a system of equations with 
B\'{e}zout number 6. Thus, whenever the values of $v_1$, $v_2$, and $v_3$
in the solutions to~\eqref{eq:threevars} are distinct, then this reflects the bound
in Lemma~\ref{le:stirling1}.

In particular, in the case $n=4$ we obtain the following 150 solutions.
\begin{description} 
\item[$k=1$, $l=1$] The six solutions for $\left(v_1,v_2,v_3\right)^T$ of the 
  system~\eqref{eq:threevars} are
\[ \left(\frac{1}{\sqrt{2}}, -\frac{1}{\sqrt{2}}, 0\right)^T, \quad
   \left(\frac{1}{20} \sqrt{110 - 30 i \sqrt{15}},
    \frac{1}{20} \sqrt{110 + 30 i \sqrt{15}},
    -\frac{1}{10}\sqrt{15} \right)^T,
\]
and the solutions obtained by permuting the first two components
of the first solution and by changing the signs and/or permuting the first 
two components in the second solution. 

For the program~\eqref{eq:invarcor} in the variables $(v_1, \ldots, v_5)^T$, 
this gives
$\binom{5}{2} \binom{2}{1} = 20$ critical positions of the form
(i.e., up to variable permutations)
\[
  \left(\frac{1}{\sqrt{2}}, -\frac{1}{\sqrt{2}}, 0, 0, 0 \right)^T \, ,
\]
20 complex solutions of the form
\[
  \left(
   -\frac{1}{20} \sqrt{110 - 30 i \sqrt{15}},
   -\frac{1}{20} \sqrt{110 + 30 i \sqrt{15}},
   \frac{1}{10} \sqrt{15}, 
   \frac{1}{10} \sqrt{15},
   \frac{1}{10} \sqrt{15}
  \right)^T \, ,
\]
and 20 complex solutions of the form
\[
  \left(
  \frac{1}{20} \sqrt{110 - 30 i \sqrt{15}},
  \frac{1}{20} \sqrt{110 + 30 i \sqrt{15}},
  -\frac{1}{10} \sqrt{15},
   -\frac{1}{10} \sqrt{15},
   -\frac{1}{10} \sqrt{15}
  \right)^T \, .
\]
\item[$k=1$, $l=2$] Here, we obtain 30 solutions of the form
\[
  \left(0,\frac{1}{2},\frac{1}{2},-\frac{1}{2},-\frac{1}{2} \right)^T \, ,
\]
30 solutions of the form
\[
  \left(  
  \frac{1}{5} \sqrt{10}, 
  \frac{1}{4}\sqrt{2} - \frac{1}{20} \sqrt{10},
  \frac{1}{4}\sqrt{2} - \frac{1}{20} \sqrt{10},
  - \frac{1}{4}\sqrt{2} - \frac{1}{20} \sqrt{10},
  -\frac{1}{4}\sqrt{2} - \frac{1}{20} \sqrt{10}
  \right)^T \, ,
\]
and 30 solutions of the form
\[
  \left(
  - \frac{1}{5} \sqrt{10},
  \frac{1}{4}\sqrt{2} + \frac{1}{20} \sqrt{10},
  \frac{1}{4}\sqrt{2} + \frac{1}{20} \sqrt{10},
  -\frac{1}{4}\sqrt{2} + \frac{1}{20} \sqrt{10},
  -\frac{1}{4}\sqrt{2} + \frac{1}{20} \sqrt{10}
  \right)^T \, .
\]
\end{description}

The global mininum is attained for the vector 
$\left(0,\frac{1}{2},\frac{1}{2},-\frac{1}{2},-\frac{1}{2} \right)^T$,
and the objective value of the global optimum is 49/80. Hence, the
radius of the smallest circumscribing cylinder for a regular simplex
in $\E^4$ with edge length $\sqrt{2}$ is 
$\sqrt{49/80} = 7 \sqrt{5} / 20 \approx 0.7826\,$.

\section*{Appendix: An error in the results of Wei{\ss}bach}

In the course of our investigations, we discovered a subtle but 
severe mistake in the paper \cite{weissbach-83-erroneous} on the explicit 
determination of the outer 
$(n{-}1)$-radius of a regular simplex in $\E^n$.
Since this error completely invalidates the proof given 
there\footnote{In a personal communication this has been confirmed
by B.~Wei{\ss}bach.},
we give a description of that flaw, including some 
computer-algebraic calculations illustrating it.

In that paper, the computation of the outer $(n{-}1)$-radius 
of a regular simplex is reduced to
the analysis of the following optimization problem.

\begin{equation}
  \label{eq:weissbach1}
  \begin{array}{rrcl}
   & \min \sum\limits_{i=1}^{n+1} u_i^4 \\
  \text{s.t.} & \sum\limits_{i=1}^{n+1} u_i^2 & = & 1 \, , \\
  & \sum\limits_{i=1}^{n+1} u_i & = & 0 \, .
  \end{array}
\end{equation}

For any local optimum $(u_1, \ldots, u_{n+1})^T$
there exist Lagrange multipliers
$\lambda_1$, $\lambda_2 \in \R$ such that 
\begin{eqnarray}
  4 u_i^3 + 2 \lambda_1 u_i + \lambda_2 \ & = & \ 0 \, , \qquad 1 \le i \le n+1 \, , \nonumber \\ 
  \sum\limits_{i=1}^{n+1} u_i^2 & = & 1 \, , \label{eq:weissbachlagrange} \\
  \sum\limits_{i=1}^{n+1} u_i & = & 0 \, . \nonumber
\end{eqnarray}
Erroneously, in \cite{weissbach-83-erroneous} it is argued 
that symmetry arguments imply that $\lambda_2 = 0$
in any solution. The following calculation in the computer algebra
system {\sc Singular} \cite{SINGULAR} shows that 
for $n=3$ this system has 26 solutions (counting multiplicity) 
over $\mathbb{C}$.

\begin{verbatim}
  ring R = 0, (u1,u2,u3,u4,la1,la2), (dp);

  ideal I =
    4*u1^3 + 2*la1*u1 + la2,
    4*u2^3 + 2*la1*u2 + la2,
    4*u3^3 + 2*la1*u3 + la2,
    4*u4^3 + 2*la1*u4 + la2,
    u1^2 + u2^2 + u3^2 + u4^2 - 1,
    u1 + u2 + u3 + u4;
  
  degree(std(I));
\end{verbatim}

This program first defines a polynomial ring in the variables 
$u_1, \ldots, u_4, \lambda_1, \lambda_2$ over a field of characteristic 
zero. We then use the {\tt degree} command to compute the dimension
and the degree of the ideal defined by our equations. The output of that
command is
\begin{verbatim}
  // codimension = 6
  // dimension   = 0
  // degree      = 26
\end{verbatim}
Hence, there are finitely many solutions (since the dimension 
of the ideal is zero), and the degree of the ideal (the sum of the 
multiplicities of the solutions) is 26.

18 of these solutions refer to the case $\lambda_2 = 0$ (and those were the 
ones computed in \cite{weissbach-83-erroneous}). Namely, if $\lambda_2 = 0$ 
then the first row of~\eqref{eq:weissbachlagrange} simplifies to
\[
  u_i (2 u_i^2 + \lambda_1) \ = \ 0 \, , \qquad 1 \le i \le n+1 \, .
\]
If we are only interested in the real solutions to this system, then
setting $\lambda_1 = -2 \lambda^2$ for some $\lambda \ge 0$ gives
\[
  u_i (u_i^2 - \lambda^2) \ = \ 0 \, , \qquad 1 \le i \le n+1 \, .
\]
Since the vector 
$(u_1, \ldots, u_{n+1})^T = (0, \ldots, 0)^T$ does not satisfy the 
second row in~\eqref{eq:weissbachlagrange}, the solutions with
$\lambda_2 = 0$ are
\begin{eqnarray*}
  u_i & = & \lambda, \qquad i \in \{i_1, \ldots, i_h\} \, ,  \\
  u_i & = & -\lambda, \qquad i \in \{i_{h+1}, \ldots, i_{2h} \} \, , \\
  u_i & = & 0 \, , \qquad i \in \{1, \ldots, n+1\} \setminus \{i_1, \ldots, i_{2h}\}
\end{eqnarray*}
for some $h \ge 1$, some set $\{i_1, \ldots, i_{2h}\}$ of pairwise different
indices, and $\lambda = (2h)^{-1/2}$.
In the case $n=3$, 
there are 12 possibilities to choose the indices and the signs 
for $|h| = 1$ and 6 possibilities to choose the indices and the signs for
$|h| = 2$, giving 18 solutions to~\eqref{eq:weissbachlagrange}.

However, there are 8 additional solutions, 
which in fact are also real! Namely, these are the solutions
\begin{eqnarray*}
  (u_1, \ldots, u_4)^T & = & 
     \frac{1}{2 \sqrt{3}} (1,-3,1,1)^T \, , \quad 
     \lambda_1 = -\frac{7}{6} \, , \quad \lambda_2 = \frac{1}{\sqrt{3}} \, , \\
  (u_1, \ldots, u_4)^T & = & 
     \frac{1}{2 \sqrt{3}} (-1,3,-1,-1)^T \, , \quad 
     \lambda_1 = -\frac{7}{6} \, , \quad \lambda_2 = - \frac{1}{\sqrt{3}} \, , \\
\end{eqnarray*}
as well as the six distinct solutions obtained from these two by permuting
the variables $u_1, \ldots, u_4$. These additional solutions invalidate
the subsequent arguments in~\cite{weissbach-83-erroneous}.

The omisssions get even worse in the higher-dimensional case. E.g.,
for $n=4$, besides the 
$\binom{5}{2} \binom{2}{1} + \binom{5}{4} \binom{4}{2} = 20 + 30 
= 50$ solutions described in \cite{weissbach-83-erroneous}, 
we obtain the following solutions:
\begin{eqnarray*}
  (u_1, \ldots, u_5)^T & = & \frac{1}{\sqrt{30}}(-2,-2,-2,3,3)^T \, , \quad 
     \lambda_1 = -\frac{7}{15} \, , \quad \, \lambda_2 = -\frac{2}{75} \sqrt{30} \, , \\
  (u_1, \ldots, u_5)^T & = & \frac{1}{\sqrt{30}}(2,2,2,-3,-3)^T \, , \quad 
     \lambda_1 = -\frac{7}{15} \, , \quad \lambda_2 = \frac{2}{75} \sqrt{30} \, , \\
  (u_1, \ldots, u_5)^T & = & \frac{1}{2 \sqrt{5}}(1,-4,1,1,1)^T \, , \quad 
     \lambda_1 = -\frac{13}{10} \, , \quad \lambda_2 = \frac{6}{25} \sqrt{5} \, , \\
  (u_1, \ldots, u_5)^T & = & \frac{1}{2 \sqrt{5}}(-1,4,-1,-1,-1)^T \, , \quad 
     \lambda_1 = -\frac{13}{10} \, , \quad \lambda_2 = -\frac{6}{25} \sqrt{5} \, ,  \\
\end{eqnarray*}
as well as those solutions obtained by permuting the variables.
Altogether, we have $10+10+5+5 = 30$ solutions with $\lambda_2 \neq 0$,
and thus a total number of 80 solutions.

Finally, we remark that the paper~\cite{weissbach-83-correct}, which computes
the outer $(n{-}1)$-radius of a regular simplex in \emph{odd} dimension $n$,
is correct (cf.\ also \cite{brandenberg-regular-2002}).

\providecommand{\bysame}{\leavevmode\hbox to3em{\hrulefill}\thinspace}
\providecommand{\MR}{\relax\ifhmode\unskip\space\fi MR }
\providecommand{\MRhref}[2]{%
  \href{http://www.ams.org/mathscinet-getitem?mr=#1}{#2}
}
\providecommand{\href}[2]{#2}

\end{document}

\bibliographystyle{amsplain}
\bibliography{bibl}
 
\end{document}